\documentclass[11pt,leqno]{amsart}
\usepackage{amsmath,amsthm,amssymb}
\usepackage{tabularx}
\newcommand{\GF}{{\mathbb F}}

\usepackage{comment}

\newcommand{\Aut}{{\rm Aut}}

\newcommand{\wt}{{\rm wt}}

\DeclareMathOperator{\Harm}{Harm}

\usepackage{color} 
\usepackage{url} 

\newtheorem{Thm}{Theorem}[section]

\newtheorem{Cor}[Thm]{Corollary}

\theoremstyle{definition}

\newtheorem{Rem}[Thm]{Remark}
\newtheorem{Ex}[Thm]{Example}

\newtheorem{Conj}[Thm]{Conjecture}

\newcommand{\NN}{\mathbb{N}}

\makeatletter
\@addtoreset{equation}{section}

\makeatother

\begin{document}

\title[A note on the Assmus--Mattson theorem]{
A note on the Assmus--Mattson theorem for some ternary codes (a resume)
}

\author{Eiichi Bannai}
\address{Faculty of Mathematics, 
Kyushu University (emeritus), Fukuoka, 819-0395 Japan}
\email{bannai@math.kyushu-u.ac.jp} 
\author{Tsuyoshi Miezaki}
\address{Faculty of Science and Engineering, 
Waseda University, 
Tokyo 169--8555, Japan
}
\email{miezaki@waseda.jp} 
\author{Hiroyuki Nakasora*}
\thanks{*Corresponding author}
\address{Faculty of Computer Science and Systems Engineering, 
Okayama Prefectural University, 
Okayama, 719-1197 Japan}
\email{nakasora@cse.oka-pu.ac.jp}

\date{}

\maketitle

\begin{abstract}
Let $C$ be a two and three-weight ternary code. 
Furthermore, we assume that $C_\ell$ are $t$-designs for all $\ell$ 
by the Assmus--Mattson theorem. 
We show that $t \leq 5$. 
As a corollary, 
we provide a new characterization of 
the (extended) ternary Golay code. 
\end{abstract}


\vspace{-0.15in}

\quad

\setcounter{section}{+0}
\section{Main results}

Let $D_{w}$ be the support design of a code $C$ for weight $w$ and
\begin{align*}
\delta(C)&:=\max\{t\in \mathbb{N}\mid \forall w, 
D_{w} \mbox{ is a } t\mbox{-design}\},\\ 
s(C)&:=\max\{t\in \mathbb{N}\mid \exists w \mbox{ s.t.~} 
D_{w} \mbox{ is a } t\mbox{-design}\}.
\end{align*}
We note that $\delta(C) \leq s(C)$.

In the present paper, we explain our main results. 
Throughout this paper, $C$ denotes a 
ternary $[n,k,d]$ code and we always assume
that a combinatorial $t$-design allows the existence of repeated blocks.
Let $C^{\perp}$ be a ternary $[n,n-k,d^{\perp}]$ dual code of $C$. 
We set $C_u:=\{c\in C\mid \wt(c)=u\}$. 
We always assume that there exists $t\in \NN$ that
satisfies the following condition: 
\begin{align}
d^{\perp}-t=\sharp\{u\mid C_u \neq \emptyset, 0<u\leq n-t\}. \label{eqn:AM}
\end{align}
This is a condition of the Assmus--Mattson theorem \cite{assmus-mattson}, 
which we call the AM-condition.
Let $D_{u}$ and $D^{\perp}_{w}$ be the support designs of $C$ and $C^{\perp}$ 
for weights $u$ and $w$, respectively. 
Then, by (\ref{eqn:AM}) and the Assmus--Mattson theorem, 
$D_u$ and $D^{\perp}_w$ are $t$-designs 
(also $s$-designs for $0<s<t$)
for any $u$ and $w$, respectively. 

Let $C$ satisfy the AM-condition.
The main results of the present paper are the following theorems. 
For a two or three-weight code,
we impose 
restrictions on $d^{\perp}$ and $t$.

\begin{Thm}\label{thm:main0}
Let $C$ be a two-weight ternary code. 
If $C$ satisfies the AM-condition, then one of the following holds: 
\begin{enumerate}
\item 
$d^\perp=5$ and $C$ is the dual of 
the ternary Golay code $[11,5,6]$ with $t=4$ or 
\item 
$d^\perp\leq 4$ and $t\leq 3$. 
\end{enumerate}
\end{Thm}

\begin{Thm}\label{thm:main2}
Let $C$ be a three-weight ternary code. 
If $C$ satisfies the AM-condition, then $d^\perp\leq 6$ and 
$t\leq 5$. 


\end{Thm}

\begin{Thm}\label{thm:main1}
Let $C$ be a three-weight ternary code, which has a weight $n$ vector. 
If $C$ satisfies the AM-condition, one of the following holds: 
\begin{enumerate}
\item 
$d^\perp=6$ and $C$ is 
the extended ternary Golay code $[12,6,6]$ with $t=5$ or 
\item 
$d^\perp\leq 5$ and $t\leq 4$. 
\end{enumerate}



\end{Thm}
It is interesting to note that 
Theorems \ref{thm:main0} (1) and \ref{thm:main1} (1) provide a 
new characterization of the (extended) ternary Golay code.


For cases in which $d^{\perp}-t=1, 2$ or $3$,
the following theorem provides a criterion for $n$ and $d$ such that 
$\delta(C^{\perp})<s(C^{\perp})$ occurs.
Let $d=d_{1}$, and $d_{2}$ and $d_{3}$ be the second and third weights of $C$, respectively.

\begin{Thm}\label{thm:main3}
Let $\alpha_{\ell}=n-d_{\ell}-(t+1)$ and $\beta_{\ell}=d_{\ell}-(t+1)$ for $\ell=1,2$ or $3$.
\begin{enumerate}
\item [{\rm (1)}]
Let $C$ satisfy the AM-condition with $d^{\perp}-t=1$.
Let $w\in \NN$ such that 
\begin{align*} 
\sum_{i+j=w} 2^{i} \binom{\alpha_{1}}{i} \cdot (-1)^{j}\binom{\beta_{1}}{j}
=0. 
\end{align*}
Then 
$D^{\perp}_{w+t+1}$ is a $(t+1)$-design if $C^{\perp}_{w+t+1}$ is non-empty. 

\item [{\rm (2)}]
Let $C$ satisfy the AM-condition with $d^{\perp}-t=2$.
Let $w\in \NN$ such that 
\begin{align*} 
\sum_{i+j=w} \left( 2^{i}  \binom{\alpha_{1}}{i} \cdot (-1)^{j} \binom{\beta_{1}}{j}
-2^{i}\binom{\alpha_{2}}{i} \cdot (-1)^{j} \binom{\beta_{2}}{j} \right) 
=0. 
\end{align*}
Then 
$D^{\perp}_{w+t+1}$ is a $(t+1)$-design if $C^{\perp}_{w+t+1}$ is non-empty. 

\item [{\rm (3)}]
Let $C$ satisfy the AM-condition with $d^{\perp}-t=3$.
Let $w\in \NN$ such that 
\begin{align*} 
\sum_{i+j=w} \biggl( 2^{i} \binom{\alpha_{1}}{i} \cdot (-1)^{j} \binom{\beta_{1}}{j}
&-\frac{d_{3}-d_{1}}{d_{3}-d_{2}} 2^{i}  \binom{\alpha_{2}}{i} \cdot (-1)^{j} \binom{\beta_{2}}{j} \\
&+\frac{d_{2}-d_{1}}{d_{3}-d_{2}} 2^{i}  \binom{\alpha_{3}}{i} \cdot (-1)^{j} \binom{\beta_{3}}{j} \biggr)
=0. 
\end{align*}
Then 
$D^{\perp}_{w+t+1}$ is a $(t+1)$-design if $C^{\perp}_{w+t+1}$ is non-empty. 

\end{enumerate}
\end{Thm}
This theorem strengthens the Assmus--Mattson theorem 
for particular cases. 
We note that parameters $n$, $d_i$, and $w$ 
that satisfy the condition in Theorem \ref{thm:main3} 
are listed on the homepage of
one of the authors \cite{miezaki}. 
In particular, 
we present the following corollary: 
\begin{Cor}\label{cor:main1}
Let $C$ satisfy the AM-condition in Theorem \ref{thm:main3}. 
For $n\leq 10$, 
in Miezaki's homepage \cite{miezaki}, 
we provide the parameters $n$, $d_i$, and $w$ 
such that $\delta(C)<s(C)$ occurs. 
\end{Cor}

Thus far, we do not possess explicit examples that fulfill Theorem~\ref{thm:main3}. 
However, for a five-weight ternary code $C$, we find an example of $\delta(C^{\perp})<s(C^{\perp})$, 
detail in Example~\ref{ex.18}.

We hold the belief that Theorem~\ref{thm:main3} holds significant importance.
In our earlier research titled ``A note on the Assmus--Mattson theorem for some binary codes," 
we introduced a binary counterpart to Theorem~\ref{thm:main3}. During that period, only a single example adhered to this theorem. 
Subsequently, we managed to identify an additional example that generates a distinct self-orthogonal design, 
characterized by a specific parameter.

For the proofs of these theorems see \cite{BMN-main} and 
for some related research, 
see 
\cite{{Bachoc},{Bachoc1},{Bannai-Koike-Shinohara-Tagami},{BM1},{BM2},{BMN},{BMY},{Delsarte},{extremal design H-M-N},{mac},{Lehmer},{Miezaki},{Miezaki2},{miezaki},{MMN},{extremal design2 M-N},{MN-tec},{dual support designs},{MN-typeIII},{mn-typeI},{Tanabe},{Venkov},{Venkov2}}. 

In Section~\ref{sec:rem}, 
we conclude the paper with remarks. 

We performed all the computer calculations in this paper with the help of 
{\sc Magma} \cite{Magma} and 
{\sc Mathematica} \cite{Mathematica}. 

\section{Concluding Remarks}\label{sec:rem}

\begin{Rem}
\begin{enumerate}

\item [(1)]


Are there examples that satisfy the 
condition of Theorem \ref{thm:main3}?

\item [(2)]

For a two-weight code, 
if we assume that $d^\perp\geq 5$ and 
\begin{align*}
W_{C} (x,y)&=x^{n}+\alpha x^{n-d_{1}}y^{d_{1}}+\beta x^{n-d_{2}}y^{d_{2}}, 
\end{align*}
then 
\[
1+\alpha+\beta=\sum_{i=0}^2\binom{n}{i}2^i. 
\]
Similarly, 
for  a three-weight code,  
if we assume that $d^\perp\geq 7$ and 
\begin{align*}
W_{C} (x,y)&=x^{n}+\alpha x^{n-d_1}y^{d_1}+\beta x^{n-d_2}y^{d_2}+\gamma x^{n-d_3}y^{d_3}, 
\end{align*}
then 
\[
1+\alpha+\beta+\gamma=\sum_{i=0}^3\binom{n}{i}2^i. 
\]
In the case $d_{3}=n$, if we assume that $d^\perp\geq 6$, 
then 
\[
1+\alpha+\beta+\gamma=3\sum_{i=0}^2\binom{n-1}{i}2^i. 
\]
In \cite{Lint}, van Lint found the solutions of the following equation for $e=2,3$:
\[
\sum_{i=0}^{e}\binom{n}{i}2^i=3^k. 
\]
Our method provides an alternative proof.


This suggests the following conjecture: 
\begin{Conj}
Let $C$ be an $\ell$-weight $[n.k.d]$ code over $\GF_{q}$ and satisfy the AM-condition.
If we assume that $d^\perp\geq 2\ell+1$ and 
\[
W_{C} (x,y)=x^{n}+\sum_{1 \leq i \leq \ell}\alpha_i x^{n-d_i}y^{d_i}, 
\]
then 
\[ 1+\alpha_1+\alpha_2+\cdots+\alpha_{\ell}=\sum_{i=0}^{\ell}\binom{n}{i}(q-1)^i =q^k. 
\]
Moreover, if $\ell \geq 4$, the codes corresponding the solutions of 
\[
\sum_{i=0}^{\ell}\binom{n}{i}(q-1)^i=q^k 
\]
do not exist. Hence, $d^\perp \leq 2\ell$ and $t \leq 2\ell-1$ for $\ell \geq 4$.

\end{Conj}

To date, we do not have a proof of this conjecture. 






\end{enumerate}

\end{Rem}

\begin{Ex}\label{ex.18}
Let $C$ be the ternary $[18,8,6]$ code which is the 158th of 160 codes in 
Database of Ternary Maximal Self-Orthogonal Codes \cite{database H-M}.
Since $\Aut(C)$ act transitively, 
the supports of any fixed weight in $C$ and $C^{\perp}$ form a $1$-design.
Furthermore $D_{9}$ and $D^{\perp}_{9}$ are $2$-designs.

We use the method of proof of Theorem \ref{thm:main3}.
The weights of $C$ are $6, 9, 12, 15, 18$, 
and $C^{\perp}$ has $d^{\perp}=4$ and $C^{\perp}_{5}$ is an empty set.
For $f \in \Harm_{2}$, 
we set
\begin{align*}
Z_{C,f}=ax^{10}y^{4} +bx^{7}y^{7}+cx^{4}y^{10}+dxy^{13}. 
\end{align*}
Then we have
\begin{align*}
Z_{C^{\perp},f}
= &a'(x+2y)^{10}(x-y)^{4} +b'(x+2y)^{7}(x-y)^{7} \\
  &c'(x+2y)^{4}(x-y)^{10} +b'(x+2y)(x-y)^{13}.
\end{align*}
Since $C^{\perp}$ has $d^{\perp}=4$ and $C^{\perp}_{5}$ is an empty set, 
the coefficients of $x^{14}$, $x^{13}y$ and $x^{11}y^{3}$  are zero in $Z_{C^{\perp},f}$. 
Then, we have $b'=0$, $c'=-3a'$ and $d'=2a'$. 
Then  we have $b=0$. Hence $D_{9}$ is a $2$-design.
Moreover 
\begin{align*}
Z_{C^{\perp},f}
= &a'\left((x+2y)^{10}(x-y)^{4} -3(x+2y)^{4}(x-y)^{10} +2(x+2y)(x-y)^{13} \right).
\end{align*}
Then the coefficients of $x^{7}y^{7}$ is zero in $Z_{C^{\perp},f}$. 
Thus $D^{\perp}_{9}$ is a $2$-design. 

\end{Ex}


\section*{Acknowledgments}

The authors would like to thank the anonymous reviewers for their 
beneficial comments on an earlier version of the manuscript. 
The second named author was supported by JSPS KAKENHI (22K03277). 
We thank Edanz (https://jp.edanz.com/ac) for editing a draft of this manuscript.


\end{document}